# UNIFORM OBSERVABILITY OF HIDDEN MARKOV MODELS AND FILTER STABILITY FOR UNSTABLE SIGNALS

By Ramon van Handel

*Princeton University*

A hidden Markov model is called observable if distinct initial laws give rise to distinct laws of the observation process. Observability implies stability of the nonlinear filter when the signal process is tight, but this need not be the case when the signal process is unstable. This paper introduces a stronger notion of uniform observability which guarantees stability of the nonlinear filter in the absence of stability assumptions on the signal. By developing certain uniform approximation properties of convolution operators, we subsequently demonstrate that the uniform observability condition is satisfied for various classes of filtering models with white-noise type observations. This includes the case of observable linear Gaussian filtering models, so that standard results on stability of the Kalman–Bucy filter are obtained as a special case.

**1. Introduction.** In a classic paper, Blackwell and Dubins [2] have obtained the following remarkably general result. Let $(Y_k)_{k\geq 0}$ be a discrete time stochastic process which takes values in a Polish space, and consider the regular conditional probabilities $\mathbf{P}((Y_k)_{k>m} \in \cdot | Y_0,\ldots,Y_m)$ and $\mathbf{Q}((Y_k)_{k>m} \in \cdot | Y_0,\ldots,Y_m)$. Then if $\mathbf{P} \sim \mathbf{Q}$, one can show that $\mathbf{P}$- and $\mathbf{Q}$-a.s.

$$\|\mathbf{P}((Y_k)_{k>m} \in \cdot | Y_0,\ldots,Y_m) - \mathbf{Q}((Y_k)_{k>m} \in \cdot | Y_0,\ldots,Y_m)\|_{\mathrm{TV}} \xrightarrow{m\to\infty} 0$$

without any further assumptions on the laws $\mathbf{P}$ and $\mathbf{Q}$. The interpretation of Blackwell and Dubins is that $\mathbf{P}$ and $\mathbf{Q}$ represent the "opinions" of two individuals about the dynamics of the time series $(Y_k)_{k\geq 0}$. When the individuals observe an initial portion of the time series $(Y_k)_{k\leq m}$, they update their opinion of the future observations $(Y_k)_{k>m}$ by Bayesian learning. The result then guarantees that the opinions of the two individuals will eventually merge, provided the individuals agree on which events can and cannot occur. A continuous time counterpart of this result was obtained by Tsukahara [27] using the prediction process of F. Knight.

The result of Blackwell and Dubins typically does not hold when **P** and **Q** are mutually singular, even when the total variation distance $\|\cdot\|_{\mathrm{TV}}$ is replaced by a weaker measure of proximity. Motivated by this problem, Diaconis and Freedman [11] investigated a special class of models with mutually singular measures for which the merging of opinions still holds in a weak sense. This has led to the investigation of various notions of *merging* of probability measures [10] which are compatible with the topology of weak convergence of probability measures. Indeed, the result of Blackwell and Dubins shows that the regular conditional probabilities $\mathbf{P}((Y_k)_{k>m} \in \cdot|(Y_k)_{k\leq m})$ and $\mathbf{Q}((Y_k)_{k>m} \in \cdot|(Y_k)_{k\leq m})$ converge toward one another, despite that neither sequence of probability measures is in fact itself convergent. Particularly when the state space is not compact, proving that two sequences of probability measures merge can be subtle (see [13], Section 11.7).

Such considerations play a central role in the present paper. Unlike the setting studied by Diaconis and Freedman, we will be content to assume the absolute continuity of our probability measures. In contrast to the problem studied by Blackwell and Dubins, however, we will consider a setting where we do not have access to the full information about the past history of the process under consideration, but we are only able to observe a subfiltration. Thus, in essence, we are interested in the merging of opinions with partial information. We will restrict ourselves to a particular aspect of this problem, the stability of the nonlinear filter, which has attracted much attention in recent years (see [9] and the references therein). As we will see, this problem can be investigated very much in the spirit of the work of Blackwell and Dubins in combination with two new ingredients: the merging of probability measures in the dual bounded-Lipschitz distance $\|\cdot\|_{\mathrm{BL}}$, as studied in the fundamental papers of Pachl [22] and Cooper and Schachermayer [7], and certain uniform approximation properties of convolution operators.

We will work chiefly in continuous time (though a general discrete time result is developed for comparison in Section 3.4). Let $(X_t, Y_t)_{t\geq 0}$ be a Markov additive process in the sense of Çinlar [6]; this means that under the probability measure $\mathbf{P}^\mu$, the processes $(X_t)_{t\geq 0}$ and $(X_t, Y_t)_{t\geq 0}$ are time-homogeneous Markov processes with initial law $(X_0, Y_0) \sim \mu \otimes \delta_0$, and that $(Y_t)_{t\geq 0}$ has conditionally independent increments given $(X_t)_{t\geq 0}$. This is the standard assumption on a hidden Markov model in continuous time, where $Y_t$ is the observed component and $X_t$ is the unobserved component. Let us now define the regular conditional probabilities $\pi_t^\mu(\cdot) = \mathbf{P}^\mu(X_t \in \cdot|(Y_s)_{s\leq t})$, that is, $\pi_t^\mu$ is the *nonlinear filter* associated to our model. We are interested in finding conditions such that $\pi_t^\mu$ and $\pi_t^\nu$ merge in an appropriate sense as $t \to \infty$ for different initial measures $\mu, \nu$.

REMARK 1.1. Using similar methods, one could also investigate the merging of the full predictive distributions $\mathbf{P}^\mu((X_r)_{r\geq t} \in \cdot|(Y_s)_{s\leq t})$. In the present paper, however, we will restrict ourselves to the study of the nonlinear filter.

Our approach has its origin in the work of Chigansky and Liptser [5], who discovered independently a corollary of the result of Blackwell and Dubins and applied it to prove that the filtered estimates of certain functions of the signal process are always stable. This idea was significantly generalized by the author in [29], where a characterization of all such functions was obtained by a functional analytic argument in the case where the signal state space is compact. In particular, it turns out that the filters $\pi_t^\mu$ and $\pi_t^\nu$ actually merge in a weak sense whenever the following *observability* condition is satisfied:

$$\mathbf{P}^\mu|_{\sigma\{(Y_r)_{r\geq 0}\}} = \mathbf{P}^\nu|_{\sigma\{(Y_r)_{r\geq 0}\}} \quad \text{implies} \quad \mu = \nu,$$

or, in other words, when distinct initial laws give rise to distinct laws of the observation process. It is tempting to conjecture that this observability criterion also leads to stability of the filter when the signal state space is not compact, as this is well known to be the case in the special case of linear Gaussian filtering models [19]. However, as the following example shows, this conjecture is not correct.

EXAMPLE 1.2. Consider a signal process $X_t$ on the state space $[1,\infty[$ defined as $X_t = X_0 e^{\lambda t}$ ($\lambda > 0$, $X_0 \geq 1$), and consider the observation process

$$Y_t = \int_0^t h(X_s)\, ds + W_t, \qquad h(x) = x^{-1}.$$

Here $W_t$ is a Wiener process independent of $X_0$. We claim that this model is observable, but that there exist $\mu \sim \nu$ such that $\pi_t^\mu$ and $\pi_t^\nu$ do not merge as $t \to \infty$.

Indeed, observability is easily demonstrated along the lines of [29], Section 5.1. To prove that $\pi_t^\mu$ and $\pi_t^\nu$ do not merge, set $f(x) = \cos(\log(x))$ and $t_n = 2\pi n/\lambda$, $n \in \mathbb{N}$. Note that $f(X_{t_n}) = f(X_0)$ for every $n \in \mathbb{N}$, so that

$$\pi_{t_n}^\rho(f) = \mathbf{E}^\rho(f(X_0)|(Y_r)_{r\leq t_n}) \xrightarrow{n\to\infty} \mathbf{E}^\rho(f(X_0)|(Y_r)_{r<\infty})$$

for any initial measure $\rho$. It thus suffices to show that $\mathbf{E}^\mu(f(X_0)|(Y_r)_{r<\infty}) \neq \mathbf{E}^\nu(f(X_0)|(Y_r)_{r<\infty})$ for some $\mu, \nu$. But by the Bayes formula

$$\mathbf{E}^\rho(f(X_0)|(Y_r)_{r<\infty})$$
$$= \frac{\int f(x)\exp(x^{-1}\int_0^\infty e^{-\lambda s}\, dY_s - 1/2 x^{-2}\int_0^\infty e^{-2\lambda s}\, ds)\rho(dx)}{\int \exp(x^{-1}\int_0^\infty e^{-\lambda s}\, dY_s - 1/2 x^{-2}\int_0^\infty e^{-2\lambda s}\, ds)\rho(dx)},$$

which is clearly not independent of $\rho$.

In the present paper we take a somewhat different point of view than in [29]. The basic idea behind our approach is easily explained. Using the Markov additive property of our model, it is not difficult to verify that

$$\mathbf{P}^\mu((Y_r - Y_t)_{r\geq t} \in \cdot|(Y_s)_{s\leq t}) = \mathbf{P}^{\pi_t^\mu}((Y_r)_{r\geq 0} \in \cdot).$$

An argument along the lines of Blackwell and Dubins applies to the left-hand side of this expression. In particular, we find that

$$\|\mathbf{P}^{\pi_t^\mu}((Y_r)_{r\geq 0} \in \cdot) - \mathbf{P}^{\pi_t^\nu}((Y_r)_{r\geq 0} \in \cdot)\|_{\mathrm{TV}} \xrightarrow{t\to\infty} 0, \qquad \mathbf{P}^\mu\text{-a.s.},$$

whenever $\mathbf{P}^\mu|_{\sigma\{(Y_r)_{r\geq 0}\}} \ll \mathbf{P}^\nu|_{\sigma\{(Y_r)_{r\geq 0}\}}$. Now suppose that we could prove that

$$\|\mathbf{P}^{\mu_n}|_{\sigma\{(Y_r)_{r\geq 0}\}} - \mathbf{P}^{\nu_n}|_{\sigma\{(Y_r)_{r\geq 0}\}}\|_{\mathrm{TV}} \xrightarrow{n\to\infty} 0 \quad \text{implies} \quad \|\mu_n - \nu_n\|_{\mathrm{BL}} \xrightarrow{n\to\infty} 0$$

for any two sequences of probability measures $\{\mu_n\}, \{\nu_n\}$. In this case, the filtering model is called *uniformly observable*, and it follows automatically that

$$\|\pi_t^\mu - \pi_t^\nu\|_{\mathrm{BL}} \xrightarrow{t\to\infty} 0, \qquad \mathbf{P}^\mu\text{-a.s.}, \qquad \text{whenever } \mathbf{P}^\mu|_{\sigma\{(Y_r)_{r\geq 0}\}} \ll \mathbf{P}^\nu|_{\sigma\{(Y_r)_{r\geq 0}\}},$$

that is, that the filters merge in the dual bounded-Lipschitz distance. This argument can be made rigorous with some care, which is done in Theorem 3.3 below.

That a filtering model may be observable but not uniformly observable is demonstrated by the counterexample above. It is easily established, however, that the two notions are identical when the signal state space is compact (Proposition 3.5), so that results of [29] follow as a special case. When the state space is not compact, proving that a filtering model is uniformly observable is more difficult. We will prove that a large class of diffusion signals with white noise type observations is uniformly observable (Section 3.4). In addition, we will show that in the linear Gaussian setting, uniform observability is equivalent to observability in the sense of linear systems theory (Section 3.3). This reproduces a well-known result on the stability of the Kalman–Bucy filter [19], which was hitherto out of reach of general stability results for nonlinear filters. The proofs of these facts rely on two key technical tools which are developed in the appendices.

The stability of nonlinear filters has been an active research topic in recent years, see, for example, [9] and the references therein. The majority of results in this direction assume that the signal process is ergodic or at least tight. Such results therefore do not allow us to prove stability of the filter when the signal process is unstable, that is, when its mass does not remain localized in a compact set. Beside specialized results for the Kalman–Bucy filter, almost all existing results in the unstable case either explicitly [4, 8, 18, 21] or implicitly [26] rely on some form of "balancing of rates" argument, where a rate of contraction must win from an opposing rate of expansion in order to give rise to stability of the filter.[1] This invariably implies that stability of the filter is only proved when the signal to noise ratio of the observations is sufficiently high. In contrast, the results in the present paper guarantee filter stability for a large class of unstable signals in a manner that is purely structural and is completely independent of the signal to noise ratio. This suggests that though one may prove filter stability by a balancing of rates argument—the latter often even leads to quantitative results on the rate of stability—this does not reflect the fundamental mechanism that causes the filter to be stable, at least in the models considered here. (The author is not aware of an example where the filter loses stability as the signal to noise ratio crosses a positive threshold.)

The remainder of this paper is organized as follows. In Section 2 we introduce the canonical hidden Markov model and the associated filtering problem. Section 3 is devoted to the statement of our main results and contains some short proofs. Longer proofs can be found in Sections 4 and 5. The appendices develop the technical tools that are used in our proofs. Appendix A establishes that certain distances between probability kernels (including the dual bounded-Lipschitz distance) are in fact measurable. Appendix B develops a general result on the merging of probability measures in the dual bounded-Lipschitz distance. This result was already obtained in a more general setting in [7, 22], but we give here a more elementary proof in the Euclidean setting. The latter is all that will be needed in our proofs, and also serves to keep the paper more self-contained. Finally, Appendix C develops a uniform approximation result for convolution operators which plays an important role in proving uniform observability for additive noise models.

**2. The hidden Markov model.** The purpose of this section is to introduce the general class of models which will be studied throughout the paper. We also introduce the filtering problem and state some fundamental regularity properties.

2.1. *Preliminaries.* Before we introduce our hidden Markov model, let us fix some notation that will be used throughout the paper.

Let $S$ be a Polish space endowed with a complete metric $d_S$. We denote by $\mathcal{B}(S)$ the Borel $\sigma$-field of $S$, and we define the spaces $B(S)$ of bounded measurable functions, $C_b(S)$ of bounded continuous functions, $U_b(S)$ of bounded uniformly continuous functions and $P(S)$ of Borel probability measures. We always endow $B(S)$, $C_b(S)$ and $U_b(S)$ with the topology of uniform convergence, and $P(S)$ with the topology of weak convergence of probability measures (recall that the space $P(S)$ is then itself Polish [23], Theorems II.6.2 and II.6.5). We denote

$$\|f\|_L = \sup_{x \neq y} \frac{|f(x) - f(y)|}{d_S(x, y)}, \qquad \|f\|_\infty = \sup_x |f(x)| \qquad \text{for all } f \in B(S),$$

and we define $\mathrm{Lip}(S) = \{f \in C_b(S) : \|f\|_\infty \leq 1 \text{ and } \|f\|_L \leq 1\}$.

---

[1]An exception is the result of [3], where filter stability is proved under the strong assumption that the observation noise has compact support. In this setting the nonlinear filter is itself compactly supported, so that this reduces essentially to the case of a compact signal state space.

Let $G \subset C_b(S)$ be uniformly bounded $\sup_{g \in G} \|g\|_\infty < \infty$, and define

$$\|\mu - \nu\|_G := \sup_{g \in G} \left| \int g \, d\mu - \int g \, d\nu \right|, \qquad \mu, \nu \in P(S).$$

Then $\|\mu - \nu\|_G$ is a pseudometric on $P(S)$, and is a metric whenever $G$ is a separating class [15], Section 3.4. We will frequently encounter the following special cases: the dual bounded-Lipschitz distance $\|\mu - \nu\|_{\mathrm{BL}} := \|\mu - \nu\|_{\mathrm{Lip}(S)}$, which metrizes the Polish space $P(S)$ [13], Theorem 11.3.3, and the total variation distance $\|\mu - \nu\|_{\mathrm{TV}} := \|\mu - \nu\|_G$ with $G = \{f \in C_b(S) : \|f\|_\infty \leq 1\}$.

As we will be interested in distances between *random* probability measures, it is important to establish that the distance $\|\mu - \nu\|_G$ is a (measurable) random variable for any pair of probability kernels $\mu, \nu$. Corollary A.2 in Appendix A establishes that this is the case whenever the family $G \subset C_b(S)$ is uniformly bounded and equicontinuous; in particular, we find that $\|\mu - \nu\|_{\mathrm{BL}}$ is measurable for any pair of probability kernels $\mu, \nu$. That the total variation distance $\|\mu - \nu\|_{\mathrm{TV}}$ between kernels is measurable is well known; this follows from the existence of a measurable version of the Radon–Nikodym derivative (see, e.g., [20], Theorem 3.1).

2.2. *Hidden Markov model.* Throughout this paper, we consider a continuous time hidden Markov model with signal state space $E$ and observation state space $\mathbb{R}^q$ (the observation dimension $q \in \mathbb{N}$ is fixed at the outset). We presume only that $E$ is Polish and we endow it with a distinguished complete metric $d$.

Let $\Omega^X = D([0, \infty[; E)$ and $\Omega^Y = D([0, \infty[; \mathbb{R}^q)$ be the spaces of $E$-valued and $\mathbb{R}^q$-valued càdlàg paths, respectively. We endow $\Omega^X$ and $\Omega^Y$ with the Skorokhod topology so that they are Polish [15], Theorem 3.5.6. We will work on the probability space $\Omega = \Omega^X \times \Omega^Y$ with its Borel $\sigma$-field $\tilde{\mathcal{F}} = \mathcal{B}(\Omega^X \times \Omega^Y)$, and we denote by $X_t : \Omega \to E$ and $Y_t : \Omega \to \mathbb{R}^q$ the coordinate projections $X_t(x, y) = x(t)$, $Y_t(x, y) = y(t)$. Furthermore, we define the natural filtrations

$$\tilde{\mathcal{F}}_t^X = \sigma\{X_s : s \leq t\}, \qquad \tilde{\mathcal{F}}_t^Y = \sigma\{Y_s : s \leq t\}, \qquad \tilde{\mathcal{F}}_t = \tilde{\mathcal{F}}_t^X \vee \tilde{\mathcal{F}}_t^Y$$

and the filtration generated by the observation increments

$$\tilde{\mathcal{G}}_t^Y = \sigma\{Y_s - Y_0 : s \leq t\}.$$

We will denote $\tilde{\mathcal{F}}^X = \bigvee_{t \geq 0} \tilde{\mathcal{F}}_t^X$, and $\tilde{\mathcal{F}}^Y$ and $\tilde{\mathcal{G}}^Y$ are defined similarly. The canonical shift $\theta_t : \Omega \to \Omega$ is defined as $\theta_t(x,y)(s) = (x(s+t), y(s+t))$.

We now proceed to impose on this canonical setup the structure of a hidden Markov model, where $Y_t$ is the observation process and $X_t$ is the signal process. Our basic assumption is that the pair $(X_t, Y_t)_{t \geq 0}$ is a time-homogeneous Markov process, whose semigroup we will denote as $T_t : B(E \times \mathbb{R}^q) \to B(E \times \mathbb{R}^q)$. We therefore presume that we are given a family $\{\mathbf{P}^\mu : \mu \in P(E)\} \subset P(\Omega)$ such that for every $\mu \in P(E)$, the pair $(X_t, Y_t)_{t \geq 0}$ is a Markov process under $\mathbf{P}^\mu$ relative to the usual augmentation [24], Section 1.4, of $\tilde{\mathcal{F}}_t$ with respect to the family $\{\mathbf{P}^\mu : \mu \in P(E)\}$, with semigroup $T_t$ and initial measure $\mu \otimes \delta_{\{0\}}$. To be precise, let us denote by $\mathcal{F}$, $\mathcal{F}^X$, $\mathcal{F}^Y$, $\mathcal{G}^Y$ the completions of $\tilde{\mathcal{F}}$, $\tilde{\mathcal{F}}^X$, $\tilde{\mathcal{F}}^Y$, $\tilde{\mathcal{G}}^Y$ and by $\mathcal{F}_t$, $\mathcal{F}_t^X$, $\mathcal{F}_t^Y$, $\mathcal{G}_t^Y$ the usual augmentations of $\tilde{\mathcal{F}}_t$, $\tilde{\mathcal{F}}_t^X$, $\tilde{\mathcal{F}}_t^Y$, $\tilde{\mathcal{G}}_t^Y$ with respect to the family $\{\mathbf{P}^\mu : \mu \in P(E)\}$. We then assume that

$$\mathbf{P}^\mu(f(X_t, Y_t)|\mathcal{F}_s) = (T_{t-s}f)(X_s, Y_s) \qquad \text{for all } f \in B(E \times \mathbb{R}^q), \mu \in P(E),$$

whenever $t \geq s \geq 0$, and that

$$\mathbf{P}^\mu(f(X_t, Y_t)) = \int (T_t f)(x, 0) \mu(dx) \qquad \text{for all } f \in B(E \times \mathbb{R}^q), \mu \in P(E).$$

Before we proceed, two remarks are in order.

REMARK 2.1. When $E$ is locally compact and $T_t$ is Feller, one can always construct the family $\mathbf{P}^\mu$ with the required properties directly from the semigroup $T_t$, for example, see [15]. As we have only assumed that $E$ is Polish, we impose the existence of the family $\mathbf{P}^\mu$ as an assumption. However, the locally compact Feller case furnishes a broad family of examples where the construction can be accomplished.

REMARK 2.2. The restriction to initial laws of the form $\mu \otimes \delta_{\{0\}}$ is in essence the requirement that the initial observation $\mathcal{F}_0^Y$ does not contain any information on the signal. The general case can be reduced to this setting, however, so there is no loss of generality in our assumptions (see the remark in [29], Section 2).

We now impose on our Markov model $(X_t, Y_t)_{t \geq 0}$ the fundamental assumption that it is a *Markov additive process* in the sense of Çinlar [6], that is, we require that the semigroup $T_t$ satisfies the following condition:

For any $f \in B(E \times \mathbb{R}^q)$, $\qquad (T_t S_y f)(x, y)$ does not depend on $y$.

Here $S_y : B(E \times \mathbb{R}^q) \to B(E \times \mathbb{R}^q)$ is defined as $(S_y f)(x, z) = f(x, z - y)$. It is not difficult to verify (see also [6]) that this assumption corresponds to the following two properties: first, the process $(X_t)_{t \geq 0}$ is a Markov process in its own right [i.e., $T_t f \in B(E)$ whenever $f \in B(E)$, where $B(E)$ is seen as a natural subspace of $B(E \times \mathbb{R}^q)$]; second, under the conditional law of $(Y_t)_{t \geq 0}$ given $\mathcal{F}^X$, the process $(Y_t)_{t \geq 0}$ has independent increments. This first property enforces the idea that there is no feedback in the system, so that the evolution of the signal is not affected by the observations. The second property enforces the idea that the observation noise is memoryless. The process $(X_t, Y_t)_{t \geq 0}$ is therefore a natural continuous time counterpart of the usual discrete time notion of a hidden Markov model, and the vast majority of continuous time filtering problems that are encountered in the literature fit in this framework (see, e.g., [30]).

2.3. *The filtering problem.* Roughly speaking, the problem of nonlinear filtering is to compute the conditional distributions $\mathbf{P}^\mu(X_t \in \cdot | \mathcal{F}_t^Y)$. As we will be dealing with convergence issues, it is essential that we choose "nice" versions of the filtered estimates. We cite the following result which provides what is needed.

LEMMA 2.3. *For every initial measure $\mu \in P(E)$, there is a probability kernel $\pi^\mu : [0, \infty[ \times \Omega \times \mathcal{B}(E) \to [0, 1]$ such that:*

1. *For every $A \in \mathcal{B}(E)$, the process $(t, \omega) \mapsto \pi^\mu(t, \omega, A)$ is the $\mathcal{F}_t^Y$-optional projection of $(t, \omega) \mapsto I_A(X_t(\omega))$.*
2. *For every $\omega \in \Omega$, the $P(E)$-valued sample path $t \mapsto \pi^\mu(t, \omega, \cdot)$ is càdlàg in the topology of $P(E)$.*

*For simplicity, we denote by $\pi_t^\mu(\cdot)$ the random measure $\omega \mapsto \pi^\mu(t, \omega, \cdot)$.*

PROOF. See [30], Proposition 1 or [17], Theorem A.3. □

As we will deal with different initial measures, the uniqueness of $\pi^\mu$ is of interest. The following result is straightforward due to the separability of $E$.

LEMMA 2.4. *The kernel $\pi^\mu$ is unique up to $\mathbf{P}^\mu |_{\mathcal{F}^Y}$-indistinguishability.*

PROOF. As $E$ is Polish, we can find a countable algebra $\{A_n\} \subset \mathcal{B}(E)$ such that $\mathcal{B}(E) = \sigma\{A_n : n \in \mathbb{N}\}$. Let $\pi^\mu$ and $\tilde{\pi}^\mu$ be two kernels that satisfy the definition of the previous lemma. To show that $\pi^\mu(t, \omega, \cdot) = \tilde{\pi}^\mu(t, \omega, \cdot)$, it suffices to show that $\pi^\mu(t, \omega, A_n) = \tilde{\pi}^\mu(t, \omega, A_n)$ for all $n$. But by the uniqueness of the optional projection up to evanescence [24], Theorem IV.5.6, we can clearly find a set $B \in \mathcal{F}^Y$ of $\mathbf{P}^\mu$-full measure such that this holds for all $t \in [0, \infty[$ and $\omega \in B$. □

**3. Main results.** The purpose of this section is to state our main results. We also give some short proofs; the remaining proofs appear in the following sections.

3.1. *Uniform observability and filter stability.* Let us begin by introducing the central result of this paper. We are interested in characterizing the *stability* of the filter, that is, the dependence of $\pi_t^\mu$ on $\mu$ as $t \to \infty$. Our general result relates this question to the following uniform notion of observability.

DEFINITION 3.1. Let $G \subset C_b(E)$ be uniformly bounded and equicontinuous. The filtering model is said to be *G-uniformly observable* if for $\{\mu_n\}, \{\nu_n\} \subset P(E)$

$$\|\mathbf{P}^{\mu_n}|_{\mathcal{F}^Y} - \mathbf{P}^{\nu_n}|_{\mathcal{F}^Y}\|_{\mathrm{TV}} \xrightarrow{n\to\infty} 0 \quad \text{implies} \quad \|\mu_n - \nu_n\|_G \xrightarrow{n\to\infty} 0.$$

When $G = \mathrm{Lip}(E)$ the model is simply called *uniformly observable*.

In [29], a model is called *observable* if $\mathbf{P}^\mu|_{\mathcal{F}^Y} = \mathbf{P}^\nu|_{\mathcal{F}^Y}$ implies $\mu = \nu$. Evidently $G$-uniform observability implies observability whenever $G$ is a separating class. However, uniform observability is strictly stronger than observability: the model is observable whenever the map $\mu \mapsto \mathbf{P}^\mu|_{\mathcal{F}^Y}$ is injective, while uniform observability requires in addition that the inverse map is uniformly continuous.[2]

REMARK 3.2. In principle one could define uniform observability in total variation by choosing $G$ to be the unit ball in $C_b(E)$ (our proofs then require some modification as this family is not equicontinuous). However, uniform observability almost never holds in this setting, as is illustrated by the following toy example.

For $\mu \in P(\mathbb{R})$, denote by $\mathbf{P}_\mu \in P(\mathbb{R})$ the law of $Y = X + \xi$ where $X \sim \mu$ and $\xi \sim N(0,1)$ are independent. Let $\mu_n = \delta_{\{1/n\}}$ and $\mu = \delta_{\{0\}}$. Then $\|\mathbf{P}_{\mu_n} - \mathbf{P}_\mu\|_{\mathrm{TV}} \to 0$ as $n \to \infty$ while $\|\mu_n - \mu\|_{\mathrm{TV}} = 2$ for all $n$. Note that this entirely reasonable model is observable and even $\mathrm{Lip}(\mathbb{R})$-uniformly observable, but uniform observability in total variation fails. Evidently we cannot obtain uniform observability in total variation when the observations are "smoothing," as is usually the case in practice, and it is therefore essential to use a smaller class $G$.

---

[2] We recall the following elementary facts. A map $f : S \to T$ between metric spaces $(S, d_S)$ and $(T, d_T)$ is called *uniformly continuous* if for every $\varepsilon > 0$, the exists a $\delta > 0$ (depending on $\varepsilon$ only) such that $d_S(x,y) < \delta$ implies $d_T(f(x), f(y)) < \varepsilon$. Equivalently, $f$ is uniformly continuous if and only if for every pair of sequences $(x_n)_{n \geq 0}$ and $(y_n)_{n \geq 0}$ such that $d_S(x_n, y_n) \to 0$, we have $d_T(f(x_n), f(y_n)) \to 0$. The proof of this fact is standard and is therefore omitted.

The following result relates the notion of uniform observability to the stability of the filter. We will prove this theorem in Section 4.

THEOREM 3.3. *Let $G \subset C_b(E)$ be uniformly bounded and equicontinuous, and suppose that the filtering model is $G$-uniformly observable. Then*

$$\|\pi_t^\mu - \pi_t^\nu\|_G \xrightarrow{t \to \infty} 0, \ \mathbf{P}^\mu\text{-}a.s., \qquad \text{whenever } \mathbf{P}^\mu|_{\mathcal{F}^Y} \ll \mathbf{P}^\nu|_{\mathcal{F}^Y}.$$

Note that in this result $G$ need not be a separating class. However, we are typically interested in the case where $G = \mathrm{Lip}(E)$. In the following subsections, we will introduce various filtering models where uniform observability can be verified.

REMARK 3.4. The condition $\mathbf{P}^\mu|_{\mathcal{F}^Y} \ll \mathbf{P}^\nu|_{\mathcal{F}^Y}$ always holds when $\mu \ll \nu$, but the latter is not necessary. It could even be the case that $\mathbf{P}^\mu|_{\mathcal{F}^Y} \sim \mathbf{P}^\nu|_{\mathcal{F}^Y}$ for every $\mu, \nu \in P(E)$, in which case the filter forgets any initial condition. The latter property is closely related to the notion of *controllability*; see [29], Section 7.

3.2. *Compact state space.* We have seen that observability in the sense of [29] is a weaker condition than uniform observability. However, in the special case that $E$ is compact and $(X, Y)$ is Feller, observability and uniform observability are equivalent. This follows directly from the general fact that any continuous bijection from a compact metric space to a metric space is a uniform homeomorphism. The proof of this fact is elementary and is given here for completeness.

PROPOSITION 3.5. *Suppose that $E$ is compact and that $(X,Y)$ is Feller. Then observability, that is, the requirement that $\mathbf{P}^\mu|_{\mathcal{F}^Y} = \mathbf{P}^\nu|_{\mathcal{F}^Y}$ implies $\mu = \nu$, already guarantees that the filtering model is uniformly observable.*

PROOF. Let $\{\mu_n\}, \{\nu_n\} \subset P(E)$ and suppose that $\|\mu_n - \nu_n\|_{\mathrm{BL}} \not\to 0$. Then we may assume, by passing to a subsequence if necessary, that $\|\mu_n - \nu_n\|_{\mathrm{BL}} \geq \varepsilon > 0$ for all $n$. As $E$ is compact, $\{\mu_n\}$ and $\{\nu_n\}$ are tight and we may assume, again passing to a subsequence if necessary, that $\|\mu_n - \mu\|_{\mathrm{BL}} \to 0$ and $\|\nu_n - \nu\|_{\mathrm{BL}} \to 0$ for some $\mu, \nu \in P(E)$. By the Feller property, we find that

$$\|\mathbf{P}^{\mu_n}|_{\mathcal{F}^Y} - \mathbf{P}^{\nu_n}|_{\mathcal{F}^Y}\|_{\mathrm{BL}} \xrightarrow{n \to \infty} \|\mathbf{P}^\mu|_{\mathcal{F}^Y} - \mathbf{P}^\nu|_{\mathcal{F}^Y}\|_{\mathrm{BL}}$$

(see [15], Theorem 4.2.5). But by the observability assumption and $\|\mu - \nu\|_{\mathrm{BL}} \geq \varepsilon$ we must have $\|\mathbf{P}^\mu|_{\mathcal{F}^Y} - \mathbf{P}^\nu|_{\mathcal{F}^Y}\|_{\mathrm{BL}} > 0$, so that $\|\mathbf{P}^{\mu_n}|_{\mathcal{F}^Y} - \mathbf{P}^{\nu_n}|_{\mathcal{F}^Y}\|_{\mathrm{TV}} \not\to 0$. By contradiction, $\|\mathbf{P}^{\mu_n}|_{\mathcal{F}^Y} - \mathbf{P}^{\nu_n}|_{\mathcal{F}^Y}\|_{\mathrm{TV}} \to 0$ must imply $\|\mu_n - \nu_n\|_{\mathrm{BL}} \to 0$. □

As a consequence, observability gives rise to stability of the filter in the sense of Theorem 3.3 when the signal state space is compact and the filtering model is Feller. Note that this result could also be obtained from the main result in [29] by using the Arzelà–Ascoli theorem (as outlined in Appendix B).

3.3. *The Kalman–Bucy filter.* Consider the hidden Markov model defined by the unique martingale problem solution to the stochastic differential equations

$$X_t = X_0 + \int_0^t AX_s \, ds + BW_t,$$

$$Y_t = \int_0^t CX_s \, ds + DV_t,$$

where $E = \mathbb{R}^d$, $A \in \mathbb{R}^{d \times d}$, $B \in \mathbb{R}^{d \times p}$, $C \in \mathbb{R}^{q \times d}$, $D \in \mathbb{R}^{q \times r}$, and $W_t$ and $V_t$ are independent $p$- and $r$-dimensional Wiener processes, respectively. We refer to this hidden Markov model as the *linear Gaussian filtering model*. When $X_0$ is Gaussian and $D$ is invertible the associated filtering problem is solved by the Kalman–Bucy filter; however, these assumptions are not required for our purposes.

We begin by stating a variant of a well-known result from linear systems theory.

LEMMA 3.6. *The following are equivalent.*

1. *The $dq \times d$-matrix*

$$O(A, C) := \begin{bmatrix} C \\ CA \\ \vdots \\ CA^{d-1} \end{bmatrix} \quad \text{has full rank.}$$

2. *There is a linear function $f : (\mathbb{R}^q)^k \to \mathbb{R}^d$ such that*

$$f\left(\int_0^{t_1} Ce^{As}x \, ds, \ldots, \int_0^{t_k} Ce^{As}x \, ds\right) = x \quad \text{for all } x \in \mathbb{R}^d$$

*for some finite number of times $t_1, \ldots, t_k \in \mathbb{R}_+$ $(k \in \mathbb{N})$.*

*When this is the case, we say that the pair $\{A, C\}$ is observable.*

PROOF. Suppose that $O(A, C)$ has full rank. We begin by noting that

$$\lim_{t \searrow 0} \frac{1}{t} \int_0^t Ce^{Ar} \, dr = C.$$

Similarly, we find that

$$\lim_{s \searrow 0} \lim_{t \searrow 0} \frac{1}{st} \left[ \int_s^{s+t} Ce^{Ar}\, dr - \int_0^t Ce^{Ar}\, dr \right] = CA.$$

Proceeding along the same lines, we can find for every $\varepsilon > 0$ a finite number of times $t_1(\varepsilon), \ldots, t_k(\varepsilon)$ and a matrix $H_\varepsilon \in \mathbb{R}^{dq \times kq}$ such that

$$H_\varepsilon \begin{bmatrix} \int_0^{t_1(\varepsilon)} Ce^{As}\, ds \\ \vdots \\ \int_0^{t_k(\varepsilon)} Ce^{As}\, ds \end{bmatrix} \xrightarrow{\varepsilon \to \infty} O(A, C).$$

But as $O(A, C)$ has full rank, the matrix on the left-hand side will have full rank for $\varepsilon$ sufficiently small, and the claim follows in one direction.

To prove the converse, note that by the Cayley–Hamilton theorem

$$\int_0^t Ce^{As}\, ds = c_0(t)C + c_1(t)CA + \cdots + c_{d-1}(t)CA^{d-1}$$

for coefficients $c_i(t)$ depending on $t$ and $A$ only. Therefore, by the existence of the function $f$, the matrix $O(A, C)$ has a left inverse and therefore has full rank. $\square$

We now obtain the following result.

PROPOSITION 3.7. *The linear Gaussian filtering model is uniformly observable if and only if $\{A, C\}$ is observable in the sense of linear systems theory.*

PROOF. We can solve the equation for $(X_t, Y_t)$ explicitly:

$$X_t = e^{At} X_0 + \int_0^t e^{A(t-s)} B\, dW_s,$$

$$Y_t = \int_0^t Ce^{As} X_0\, ds + \int_0^t \int_0^s Ce^{A(s-r)} B\, dW_r\, ds + DV_t.$$

Suppose first that $\{A, C\}$ is not observable. Then there exists $v \in \mathbb{R}^d$ such that

$$\int_0^t Ce^{As} v\, ds = 0 \qquad \text{for all } t \geq 0.$$

When this is the case, it is easily seen that for any initial law $\mu \in P(\mathbb{R}^d)$, the initial law $\mu * \delta_v$ gives rise to the same law of the observations as does $\mu$. Therefore the model is certainly not uniformly observable.

Conversely, suppose that the pair $\{A, C\}$ is observable. Let $t_1, \ldots, t_k$ and $f : (\mathbb{R}^q)^k \to \mathbb{R}^d$ be as in Lemma 3.6. Then we can write

$$(Y_{t_1}, \ldots, Y_{t_k}) = \left( \int_0^{t_1} Ce^{As} X_0 \, ds, \ldots, \int_0^{t_k} Ce^{As} X_0 \, ds \right) + \xi,$$

where $\xi$ is a $kq$-dimensional Gaussian random variable. In particular, the characteristic function of $\xi$ vanishes nowhere. By Proposition C.2 and the fact that $f$ is Lipschitz continuous (as it is linear), it is easily established that

$$\|\mathbf{P}^{\mu_n}|_{\mathcal{F}^Y} - \mathbf{P}^{\nu_n}|_{\mathcal{F}^Y}\|_{\mathrm{TV}} \xrightarrow{n \to \infty} 0 \quad \text{implies} \quad \|\mu_n - \nu_n\|_{\mathrm{BL}} \xrightarrow{n \to \infty} 0.$$

This completes the proof of uniform observability. □

As a corollary, it follows from Theorem 3.3 that if $\{A, C\}$ is observable, then $\|\pi_t^\mu - \pi_t^\nu\|_{\mathrm{BL}} \to 0$ $\mathbf{P}^\mu$-a.s. as $t \to \infty$ whenever $\mathbf{P}^\mu|_{\mathcal{F}^Y} \ll \mathbf{P}^\nu|_{\mathcal{F}^Y}$. This result is essentially known, see, for example, [19], Section 2. However, previous proofs rely crucially on the fact that the solution to the filtering problem can be explicitly expressed in terms of the Kalman–Bucy filtering equations, which are amenable to explicit analysis. In contrast, the Kalman–Bucy filter (in the case of unstable signals) has hitherto been out of reach of results on filter stability which also apply to nonlinear filtering models. The present approach is therefore of significant interest, as it allows us to infer stability of the filter directly from the general Theorem 3.3.

REMARK 3.8. The present result differs somewhat from previous stability results for the Kalman–Bucy filter. It is customary to assume *controllability* in addition to observability, which is replaced in our setting by the absolute continuity requirement $\mathbf{P}^\mu|_{\mathcal{F}^Y} \ll \mathbf{P}^\nu|_{\mathcal{F}^Y}$. It is not difficult to verify that if the signal is controllable and $D$ is invertible, then $\mathbf{P}^\mu|_{\mathcal{F}^Y} \sim \mathbf{P}^\nu|_{\mathcal{F}^Y}$ for every $\mu, \nu \in P(\mathbb{R}^d)$, so that our result is in fact more general in this sense. On the other hand, the assumptions in [19], Section 2, are weaker than the observability assumption; in particular, *detectability* suffices (at least when $D$ is invertible; see also [29], Appendix A). It would be interesting to obtain a generalization of the latter notion to general hidden Markov models, for example, by combining Theorem 3.3 with the results in [28].

3.4. *Diffusion signals.* The verification of uniform observability for the linear Gaussian filtering model was simplified significantly by the fact that the stochastic differential equations which define the model can be solved explicitly. In the present subsection we will verify uniform observability for a class of *nonlinear* filtering models, where we do not have this luxury. Consequently the conditions for uniform observability will be more stringent than in the previous section; in particular, we will recover Proposition 3.7

14  R. VAN HANDELas a special case in the setting where $C$ is invertible (in which case $\{A, C\}$ is automatically observable).

Let $E = \mathbb{R}^q$ (i.e., we assume that the signal and observation state space dimensions coincide). We consider a hidden Markov model of the form

$$X_t = X_0 + \int_0^t b(X_s)\,ds + \int_0^t \sigma(X_s)\,dW_s,$$

$$Y_t = \int_0^t h(X_s)\,ds + DV_t,$$

where $W_t$ and $V_t$ are independent $p$- and $r$-dimensional Wiener processes, respectively, and $D \in \mathbb{R}^{q \times r}$, $b : \mathbb{R}^q \to \mathbb{R}^q$, $\sigma : \mathbb{R}^q \to \mathbb{R}^{q \times p}$, $h : \mathbb{R}^q \to \mathbb{R}^q$. In addition, we assume that the following conditions hold:

1. $b$ is globally Lipschitz continuous;
2. $\sigma$ is globally Lipschitz continuous and bounded;
3. $h(x) = Cx + h_0(x)$, where $C$ is an invertible matrix and $\|C^{-1}h_0\|_L < 1$.

Note that under these conditions it is well known that the martingale problem for $(X, Y)$ has a unique solution, so that our model is well defined.

The proof of the following result can be found in Section 5.

THEOREM 3.9. *The filtering model in this section is uniformly observable.*

The required form of the observation function $h$ may seem a little odd; however, the proof of Theorem 3.9 shows that this is a natural choice. To gain a little more insight into this condition, we prove the following lemma.

LEMMA 3.10. *Any function $h(x) = Cx + h_0(x)$, where $C$ is invertible and $\|C^{-1}h_0\|_L < 1$, is bi-Lipschitz, that is, there exist $0 < m < M < \infty$ such that*

$$m\|x - y\| \leq \|h(x) - h(y)\| \leq M\|x - y\| \qquad \text{for all } x, y \in \mathbb{R}^q.$$

*Conversely, if $q = 1$ and $h$ is a bi-Lipschitz function, then $h(x) = Cx + h_0(x)$ for some $0 < C < \infty$ and Lipschitz function $h_0$ with $\|C^{-1}h_0\|_L < 1$.*

PROOF. Suppose that $h(x) = Cx + h_0(x)$, where $C$ is an invertible matrix and $\|C^{-1}h_0\|_L < 1$. Clearly $M := \|h\|_L < \infty$. Moreover, we can estimate

$$\|x - y\| \leq \|C^{-1}h(x) - C^{-1}h(y)\| + \|C^{-1}h_0(x) - C^{-1}h_0(y)\|$$
$$\leq \|C^{-1}\|\|h(x) - h(y)\| + \|C^{-1}h_0\|_L\|x - y\|.$$

As $\|C^{-1}h_0\|_L < 1$, we may set $m := (1 - \|C^{-1}h_0\|_L)/\|C^{-1}\|$.

Conversely, let $q = 1$ and suppose that $h$ is bi-Lipschitz with constants $m < M$. Then in particular $h: \mathbb{R} \to \mathbb{R}$ is a continuous bijection, so that it is either strictly increasing or strictly decreasing. Define $C := (M+m)/2$ if $h$ is increasing and $C := -(M+m)/2$ if $h$ is decreasing. Then for any $x > y$, we evidently have
$$(1-\varepsilon)(x-y) \leq C^{-1}h(x) - C^{-1}h(y) \leq (1+\varepsilon)(x-y),$$
where $\varepsilon := (M-m)/(M+m)$. In particular,
$$\frac{|C^{-1}h(x) - x - (C^{-1}h(y) - y)|}{|x-y|} \leq \varepsilon \qquad \text{for all } x > y.$$
This estimate consequently holds for all $x, y \in \mathbb{R}$ by symmetry. The result now follows by noting that $h_0(x) := h(x) - Cx$ satisfies $\|C^{-1}h_0\|_L \leq \varepsilon < 1$.
□

Using Lemma 3.10 we find that when the signal state space is the real line, the filtering model of the present section is uniformly observable whenever the observation function $h$ is a Lipschitz bijection with Lipschitz inverse (i.e., bi-Lipschitz). In higher dimensions the condition $h(x) = Cx + h_0(x)$ is stronger than the bi-Lipschitz condition, and enforces the idea that $h(x)$ cannot be "too nonlinear."

Intuitively, one might well expect that for additive noise observation models with a strongly invertible observation function $h$, the filter would be stable under only mild conditions on the signal process. This is certainly the spirit of Theorem 3.9, but the requirement on $h$ and the assumptions on the signal (i.e., that it is a diffusion) are somewhat stronger than one might expect to be necessary. Following the approach used in the proof of Theorem 3.9, the author did not succeed in weakening the requirements of that result. For comparison, however, let us briefly discuss a related problem in discrete time where a very general result may be obtained.

Let $E = \mathbb{R}^q$, and let $P: \mathbb{R}^q \times \mathcal{B}(\mathbb{R}^q) \to [0,1]$ be a given transition probability kernel. On the sequence space $E^{\mathbb{Z}_+} \times F^{\mathbb{Z}_+}$ with the canonical coordinate projections $X_n(x, y) = x(n)$, $Y_n(x, y) = y(n)$, we define the family of probability measures $\mathbf{P}^\mu$, $\mu \in P(\mathbb{R}^q)$ such that $(X_n)_{n \geq 0}$ is a Markov chain with initial measure $X_0 \sim \mu$ and transition probability $P$, and such that $Y_n = h(X_n) + \xi_n$ for every $n \geq 0$ where $\xi_n$ is an i.i.d. sequence independent of $(X_n)_{n \geq 0}$. We now define for every $\mu \in P(\mathbb{R}^q)$ the regular conditional probabilities
$$\pi_n^\mu(\cdot) := \mathbf{P}^\mu(X_{n+1} \in \cdot | Y_0, \ldots, Y_n), \qquad n \geq 0.$$
In other words, $\pi_n^\mu$ is the *one step predictor* of the signal given the observations.

In the present setting, the following result holds without further assumptions.

PROPOSITION 3.11. *Suppose that the following hold:*

1. *$h$ possesses a uniformly continuous inverse; and*
2. *the characteristic function of $\xi_0$ vanishes nowhere.*

*Then $\|\pi_n^\mu - \pi_n^\nu\|_{\mathrm{BL}} \xrightarrow{n\to\infty} 0$, $\mathbf{P}^\mu$-a.s. whenever $\mathbf{P}^\mu|_{\sigma\{(Y_k)_{k\geq 0}\}} \ll \mathbf{P}^\nu|_{\sigma\{(Y_k)_{k\geq 0}\}}$.*

PROOF. Let $\mu, \nu \in P(\mathbb{R}^q)$ satisfy $\mathbf{P}^\mu|_{\sigma\{(Y_k)_{k\geq 0}\}} \ll \mathbf{P}^\nu|_{\sigma\{(Y_k)_{k\geq 0}\}}$, and let $\xi \in P(\mathbb{R}^q)$ be the law of $\xi_0$. It is easily verified that for any $\rho \in P(\mathbb{R}^q)$

$$\mathbf{P}^\rho(Y_{n+1} \in \cdot | Y_0, \ldots, Y_n) = \pi_n^\rho h^{-1} * \xi.$$

The classical result of Blackwell and Dubins [2], Section 2, shows that

$$\|\pi_n^\mu h^{-1} * \xi - \pi_n^\nu h^{-1} * \xi\|_{\mathrm{TV}} \xrightarrow{n\to\infty} 0, \qquad \mathbf{P}^\mu\text{-a.s.}$$

We therefore obtain by Proposition C.2

$$\|\pi_n^\mu h^{-1} - \pi_n^\nu h^{-1}\|_{\mathrm{BL}} \xrightarrow{n\to\infty} 0, \qquad \mathbf{P}^\mu\text{-a.s.}$$

As the bounded-Lipschitz functions are uniformly dense in $U_b(\mathbb{R}^q)$ [12], Lemma 8,

$$|\pi_n^\mu(f \circ h) - \pi_n^\nu(f \circ h)| \xrightarrow{n\to\infty} 0 \quad \text{for all } f \in U_b(\mathbb{R}^q),\ \mathbf{P}^\mu\text{-a.s.},$$

where the $\mathbf{P}^\mu$-exceptional set does not depend on $f$. But $h$ has a uniformly continuous inverse, so any function in $U_b(\mathbb{R}^q)$ can be written as $f \circ h$ for some $f \in U_b(\mathbb{R}^q)$. The result now follows from Corollary B.4.  □

It should be noted, in particular, that this result places no conditions whatsoever on the signal process $X_n$ except the Markov property. However, this result is a statement about the one step predictor and not about the filter. In continuous time, one can obtain filtered estimates at time $t$ by taking the limit of predictive estimates over the time interval $[t, t+\delta]$ as $\delta \searrow 0$. The chief difficulty in the proof of Theorem 3.9 is to show that the limits as $\delta \searrow 0$ and $t \to \infty$ can be interchanged.

**4. Proof of Theorem 3.3.** In the following, we denote by $F^Y$ the family

$$F^Y = \mathrm{span}\{f_1(Y_{t_1} - Y_0) \cdots f_k(Y_{t_k} - Y_0):$$
$$f_1, \ldots, f_k \in B(\mathbb{R}^q), t_1, \ldots, t_k \in [0, \infty[, k \in \mathbb{N}\}$$

of $\tilde{\mathcal{G}}^Y$-measurable cylindrical random variables. Before we turn to the proof of Theorem 3.3, we introduce two elementary lemmas.

LEMMA 4.1. *There is a countable $H^Y \subset F^Y$, $\sup_{h \in H^Y} \|h\|_\infty \leq 1$ so that*

$$\|\mathbf{P}^\mu|_{\mathcal{F}^Y} - \mathbf{P}^\nu|_{\mathcal{F}^Y}\|_{\mathrm{TV}} = \sup_{h \in H^Y} \left| \int h\, d\mathbf{P}^\mu - \int h\, d\mathbf{P}^\nu \right| \qquad \text{for all } \mu, \nu \in P(E).$$

PROOF. Note that for any $\mu \in P(E)$, the $\sigma$-fields $\mathcal{F}^Y$ and $\tilde{\mathcal{G}}^Y$ coincide $\mathbf{P}^\mu$-a.s. By [15], Proposition 3.7.1, we have $\tilde{\mathcal{G}}^Y = \bigvee_{k \geq 0} \tilde{\mathcal{G}}^{Y,k}$ where

$$\tilde{\mathcal{G}}^{Y,k} = \sigma\{Y_{2^{-k}\ell} - Y_0 : \ell = 1, \ldots, 4^k\}.$$

Choose a countable dense set $\{x_p\} \subset \mathbb{R}^q$, and consider the countable collection of open balls $B_{p,m} = \{x \in \mathbb{R}^q : |x - x_p| < 1/m\}$. Then $\tilde{\mathcal{G}}^{Y,k} = \bigvee_{n \geq 0} \tilde{\mathcal{G}}^{Y,k,n}$ with

$$\tilde{\mathcal{G}}^{Y,k,n} = \sigma\{Y_{2^{-k}\ell} - Y_0 \in B_{p_\ell, m_\ell} : p_\ell, m_\ell = 1, \ldots, n, \ell = 1, \ldots, 4^k\}.$$

Now note that every $\tilde{\mathcal{G}}^{Y,k,n}$ consists of a finite number of sets in $\tilde{\mathcal{G}}^Y$, and for every $A \in \tilde{\mathcal{G}}^{Y,k,n}$ the indicator function $I_A \in F^Y$. But $\tilde{\mathcal{G}}^{Y,k,n} \nearrow \tilde{\mathcal{G}}^{Y,k}$ as $n \to \infty$ and $\tilde{\mathcal{G}}^{Y,k} \nearrow \tilde{\mathcal{G}}^Y$ as $k \to \infty$, so that we can evidently estimate

$$\|\mathbf{P}^\mu|_{\mathcal{F}^Y} - \mathbf{P}^\nu|_{\mathcal{F}^Y}\|_{\mathrm{TV}} = \|\mathbf{P}^\mu|_{\tilde{\mathcal{G}}^Y} - \mathbf{P}^\nu|_{\tilde{\mathcal{G}}^Y}\|_{\mathrm{TV}}$$

$$= 2 \lim_{k \to \infty} \lim_{n \to \infty} \max_{A \in \tilde{\mathcal{G}}^{Y,k,n}} |\mathbf{P}^\mu(A) - \mathbf{P}^\nu(A)|$$

$$\leq \sup_{h \in H^Y} \left| \int h\, d\mathbf{P}^\mu - \int h\, d\mathbf{P}^\nu \right|,$$

where we have defined the countable family $H^Y \subset F^Y$ as

$$H^Y = \bigcup_{k,n \in \mathbb{N}} \{I_A - I_{A^c} : A \in \tilde{\mathcal{G}}^{Y,k,n}\}.$$

On the other hand, the reverse inequality is immediate. □

We will also need the following.

LEMMA 4.2. *For any $\xi \in F^Y$, $\mu \in P(E)$ and $t \in [0, \infty[$, we have*

$$\mathbf{E}^\mu(\xi \circ \theta_t | \mathcal{F}_t^Y) = \mathbf{E}^{\pi_t^\mu}(\xi), \qquad \mathbf{P}^\mu\text{-a.s.}$$

PROOF. By the Markov additive property of our model,

$$\mathbf{E}^\mu(f(Y_{s+t} - Y_t)|\mathcal{F}_t) = \mathbf{E}^\mu(f(Y_{s+t} - y)|\mathcal{F}_t)|_{y=Y_t}$$
$$= T_s S_y f(X_t, Y_t)|_{y=Y_t} = T_s f(X_t, 0) = \mathbf{E}^{\delta_{X_t}}(f(Y_s - Y_0)).$$

Along the same lines $\mathbf{E}^\mu(\xi \circ \theta_t | \mathcal{F}_t) = \mathbf{E}^{\delta_{X_t}}(\xi)$ for any $\xi \in F^Y$. Therefore

$$\mathbf{E}^\mu(\xi \circ \theta_t | \mathcal{F}_t^Y) = \mathbf{E}^\mu(\mathbf{E}^{\delta_{X_t}}(\xi)|\mathcal{F}_t^Y) = \int \mathbf{E}^{\delta_x}(\xi) \pi_t^\mu(dx) = \mathbf{E}^{\pi_t^\mu}(\xi)$$

by the tower property of the conditional expectation, and the result follows. □

We now turn to the proof of Theorem 3.3. We assume throughout the proof that $\mathbf{P}^\mu|_{\mathcal{F}^Y} \ll \mathbf{P}^\nu|_{\mathcal{F}^Y}$, so that in particular both $\pi_t^\mu$ and $\pi_t^\nu$ are defined uniquely $\mathbf{P}^\mu$-a.s.

PROOF OF THEOREM 3.3. Let $\xi \in H^Y$; then we have

$$|\mathbf{E}^{\pi_t^\mu}(\xi) - \mathbf{E}^{\pi_t^\nu}(\xi)| = \frac{|\mathbf{E}^\nu((\Lambda - \mathbf{E}^\nu(\Lambda|\mathcal{F}_t^Y))\xi \circ \theta_t|\mathcal{F}_t^Y)|}{\mathbf{E}^\nu(\Lambda|\mathcal{F}_t^Y)}, \qquad \mathbf{P}^\mu\text{-a.s.}$$

by Lemma 4.2 and the Bayes formula, where $\Lambda := d\mathbf{P}^\mu|_{\mathcal{F}^Y}/d\mathbf{P}^\nu|_{\mathcal{F}^Y}$. But $H^Y$ is countable and $\sup_{\xi \in H^Y} \|\xi\|_\infty \leq 1$, so we can evidently estimate

$$\sup_{\xi \in H^Y} |\mathbf{E}^{\pi_t^\mu}(\xi) - \mathbf{E}^{\pi_t^\nu}(\xi)| \leq \frac{\mathbf{E}^\nu(M_t|\mathcal{F}_t^Y)}{\mathbf{E}^\nu(\Lambda|\mathcal{F}_t^Y)} \qquad \text{for all } t \in \mathbb{Q}_+, \ \mathbf{P}^\mu\text{-a.s.,}$$

where $M_t := |\Lambda - \mathbf{E}^\nu(\Lambda|\mathcal{F}_t^Y)|$. It should be noted that as $\mathbf{P}^\mu|_{\mathcal{F}^Y} \ll \mathbf{P}^\nu|_{\mathcal{F}^Y}$, all the preceding quantities are $\mathbf{P}^\mu|_{\mathcal{F}^Y}$-a.s. uniquely defined and we have implicitly taken only countable intersections of sets of full measure (as $H^Y$ and $\mathbb{Q}_+$ are countable).

We now claim that $\mathbf{E}^\nu(M_t|\mathcal{F}_t^Y) \to 0$ $\mathbf{P}^\nu$-a.s. as $t \to \infty$ along the rationals $t \in \mathbb{Q}_+$. To see this, define $M_t^k := |\Lambda I_{\Lambda \leq k} - \mathbf{E}^\nu(\Lambda I_{\Lambda \leq k}|\mathcal{F}_t^Y)|$, and note that

$$M_t \leq M_t^k + \Lambda I_{\Lambda > k} + \mathbf{E}^\nu(\Lambda I_{\Lambda > k}|\mathcal{F}_t^Y) \qquad \text{for all } t \in \mathbb{Q}_+, k \in \mathbb{N}, \ \mathbf{P}^\nu\text{-a.s.}$$

Therefore we obtain, using that trivially $2\Lambda I_{\Lambda > k} \to 0$ as $k \to \infty$ $\mathbf{P}^\nu$-a.s.,

$$\limsup_{t \to \infty, t \in \mathbb{Q}_+} \mathbf{E}^\nu(M_t|\mathcal{F}_t^Y) \leq \limsup_{k \to \infty} \limsup_{t \to \infty, t \in \mathbb{Q}_+} \mathbf{E}^\nu(M_t^k|\mathcal{F}_t^Y), \qquad \mathbf{P}^\nu\text{-a.s.}$$

But as by construction $\mathbf{P}^\nu$-a.s. $M_t^k \leq k$ for all $t \in \mathbb{Q}_+$, $k \in \mathbb{N}$ and as $\mathbf{P}^\nu$-a.s. $\limsup_{t \to \infty, t \in \mathbb{Q}_+} M_t^k = 0$ for all $k \in \mathbb{N}$ by martingale convergence, we have

$$\limsup_{t \to \infty, t \in \mathbb{Q}_+} \mathbf{E}^\nu(M_t^k|\mathcal{F}_t^Y) \leq \limsup_{n \to \infty} \limsup_{t \to \infty, t \in \mathbb{Q}_+} \mathbf{E}^\nu\left(\sup_{\substack{s \geq n \\ s \in \mathbb{Q}_+}} M_s^k \Big| \mathcal{F}_t^Y\right) = 0, \qquad \mathbf{P}^\nu\text{-a.s.}$$

As by construction $\Lambda > 0$ $\mathbf{P}^\mu$-a.s., we have evidently established that

$$\sup_{\xi \in H^Y} |\mathbf{E}^{\pi_t^\mu}(\xi) - \mathbf{E}^{\pi_t^\nu}(\xi)| \xrightarrow[t \in \mathbb{Q}_+]{t \to \infty} 0, \qquad \mathbf{P}^\mu\text{-a.s.}$$

Denote by $\Omega_0 \subset \Omega$ a set of $\mathbf{P}^\mu$-full measure on which this convergence holds. Then

$$\|\mathbf{P}^{\pi^\mu(t_k, \omega, \cdot)}|_{\mathcal{F}^Y} - \mathbf{P}^{\pi^\nu(t_k, \omega, \cdot)}|_{\mathcal{F}^Y}\|_{\text{TV}} \xrightarrow{k \to \infty} 0$$

for every $\omega \in \Omega_0$ and every subsequence $\{t_k\} \subset \mathbb{Q}_+$ such that $t_k \nearrow \infty$. As the model is presumed to be $G$-uniformly observable, this implies that

$$\|\pi^\mu(t_k, \omega, \cdot) - \pi^\nu(t_k, \omega, \cdot)\|_G \xrightarrow{k \to \infty} 0$$

for every $\omega \in \Omega_0$ and every subsequence $\{t_k\} \subset \mathbb{Q}_+$ such that $t_k \nearrow \infty$. But as $G$ is uniformly bounded and equicontinuous and as $\pi^\mu(t,\omega,\cdot)$ and $\pi^\nu(t,\omega,\cdot)$ are càdlàg in the topology of $P(E)$ by Lemma 2.3, it follows from [23], Theorem II.6.8, that $t \mapsto \|\pi_t^\mu - \pi_t^\nu\|_G$ is càdlàg. We therefore obtain

$$\|\pi^\mu(t,\omega,\cdot) - \pi^\nu(t,\omega,\cdot)\|_G \xrightarrow{t\to\infty} 0 \qquad \text{for all } \omega \in \Omega_0.$$

This completes the proof. □

**5. Proof of Theorem 3.9.** In the proof of Theorem 3.9 we will make essential use of the flow generated by the deterministic part of the signal process: define $\eta_t(x)$, for every $x \in \mathbb{R}^q$, as the solution of the ordinary differential equation

$$\eta_t(x) = x + \int_0^t b(\eta_s(x))\,ds.$$

Existence and uniqueness follows from the global Lipschitz property of $b$.

The special form of $h$ is essential, as it allows us to establish the following.

LEMMA 5.1. *Let $h(x) = Cx + h_0(x)$, where $C$ is an invertible matrix and $\|C^{-1}h_0\|_L < 1$. Then there exist constants $\varepsilon_0 > 0$ and $m, M > 0$ such that*

$$m\|x-y\| \le \left\|\frac{1}{\varepsilon}\int_0^\varepsilon h(\eta_s(x))\,ds - \frac{1}{\varepsilon}\int_0^\varepsilon h(\eta_s(y))\,ds\right\| \le M\|x-y\|$$

*for every $\varepsilon < \varepsilon_0$ and $x,y \in \mathbb{R}^q$.*

PROOF. Let us define

$$H_\varepsilon(x) := \frac{1}{\varepsilon}\int_0^\varepsilon h(\eta_s(x))\,ds$$

and note that we can write

$$C^{-1}H_\varepsilon(x) = \frac{1}{\varepsilon}\int_0^\varepsilon \eta_s(x)\,ds + \frac{1}{\varepsilon}\int_0^\varepsilon C^{-1}h_0(\eta_s(x))\,ds.$$

We now estimate as follows.

$$\|x-y\| \le \|C^{-1}H_\varepsilon(x) - C^{-1}H_\varepsilon(y)\| + \|C^{-1}H_\varepsilon(x) - C^{-1}H_\varepsilon(y) - (x-y)\|$$

$$\le \|C^{-1}H_\varepsilon(x) - C^{-1}H_\varepsilon(y)\| + \frac{1}{\varepsilon}\int_0^\varepsilon \|\eta_s(x) - \eta_s(y) - (x-y)\|\,ds$$

$$+ \frac{1}{\varepsilon}\int_0^\varepsilon \|C^{-1}h_0(\eta_s(x)) - C^{-1}h_0(\eta_s(y))\|\,ds$$

$$\le \|C^{-1}\|\|H_\varepsilon(x) - H_\varepsilon(y)\| + \frac{1}{\varepsilon}\int_0^\varepsilon \int_0^s \|b(\eta_r(x)) - b(\eta_r(y))\|\,dr\,ds$$

$$+ \|C^{-1}h_0\|_L \frac{1}{\varepsilon} \int_0^\varepsilon \|\eta_s(x) - \eta_s(y)\| \, ds$$

$$\leq \|C^{-1}\| \|H_\varepsilon(x) - H_\varepsilon(y)\|$$
$$+ (\|C^{-1}h_0\|_L + \|b\|_L \varepsilon/2) \sup_{s<\varepsilon} \|\eta_s(x) - \eta_s(y)\|.$$

But note that

$$\|\eta_s(x) - \eta_s(y)\| \leq \|x - y\| + \|b\|_L \int_0^s \|\eta_r(x) - \eta_r(y)\| \, dr,$$

so that by Gronwall's lemma

$$\sup_{s<\varepsilon} \|\eta_s(x) - \eta_s(y)\| \leq e^{\|b\|_L \varepsilon} \|x - y\|.$$

We therefore find that for all $x, y \in \mathbb{R}^q$ and $\varepsilon > 0$

$$\frac{1 - \|C^{-1}h_0\|_L e^{\|b\|_L \varepsilon} - \|b\|_L \varepsilon e^{\|b\|_L \varepsilon}/2}{\|C^{-1}\|} \|x - y\| \leq \|H_\varepsilon(x) - H_\varepsilon(y)\|.$$

But evidently

$$\frac{1 - \|C^{-1}h_0\|_L e^{\|b\|_L \varepsilon} - \|b\|_L \varepsilon e^{\|b\|_L \varepsilon}/2}{\|C^{-1}\|} \nearrow \frac{1 - \|C^{-1}h_0\|_L}{\|C^{-1}\|} > 0 \quad \text{as } \varepsilon \searrow 0.$$

This establishes the lower bound. For the upper bound, note that

$$\|H_\varepsilon(x) - H_\varepsilon(y)\| \leq \frac{1}{\varepsilon} \int_0^\varepsilon \|h(\eta_s(x)) - h(\eta_s(y))\| \, ds$$
$$\leq \|h\|_L \sup_{s \leq \varepsilon} \|\eta_s(x) - \eta_s(y)\| \leq \|h\|_L e^{\|b\|_L \varepsilon} \|x - y\|.$$

The proof is complete. □

The following lemma is used to reduce the proof of Theorem 3.9 to the study of the deterministic part $\eta_t(x)$, rather than working with the fully stochastic signal $X_t$. It is here that the boundedness of the diffusion coefficient $\sigma$ is used.

LEMMA 5.2. *Provided that $\sigma$ is bounded, we have*

$$\sup_{s \leq t} \sup_{\mu \in P(\mathbb{R}^q)} \mathbf{E}^\mu(\|X_s - \eta_s(X_0)\|) \xrightarrow{t \to 0} 0.$$

PROOF. For every $x \in \mathbb{R}^q$, let $\xi_t(x)$ be the solution of

$$\xi_t(x) = x + \int_0^t b(\xi_s(x)) \, ds + \int_0^t \sigma(\xi_s(x)) \, dW_s.$$

By the global Lipschitz property of the coefficients, the solution is uniquely defined and is square integrable for every $x \in \mathbb{R}^q$. We therefore obtain using Itô's rule

$$\mathbf{E}(\|\xi_t(x) - \eta_t(x)\|^2)$$
$$= \mathbf{E}\left[\int_0^t \{2\langle \xi_s(x) - \eta_s(x), b(\xi_s(x)) - b(\eta_s(x))\rangle + a(\xi_s(x))\} \, ds\right],$$

where $a(x) = \text{Tr}[\sigma(x)^*\sigma(x)]$. Note that as we have assumed that $\sigma$ is uniformly bounded, $a(x)$ is also uniformly bounded $a(x) \leq K < \infty$. Therefore

$$\mathbf{E}(\|\xi_t(x) - \eta_t(x)\|^2) \leq Kt + 2\|b\|_L \int_0^t \mathbf{E}(\|\xi_s(x) - \eta_s(x)\|^2) \, ds.$$

By Gronwall's lemma, we obtain for every $T < \infty$ and $x \in \mathbb{R}^q$

$$\sup_{t \leq T} \mathbf{E}(\|\xi_t(x) - \eta_t(x)\|^2) \leq KTe^{2\|b\|_L T}.$$

By Jensen's inequality, we find that for every $T < \infty$

$$\sup_{t \leq T} \sup_{x \in \mathbb{R}^q} \mathbf{E}(\|\xi_t(x) - \eta_t(x)\|) \leq e^{\|b\|_L T}\sqrt{KT}.$$

It remains to note that

$$\sup_{\mu \in P(\mathbb{R}^q)} \mathbf{E}^\mu(\|X_t - \eta_t(X_0)\|) = \sup_{\mu \in P(\mathbb{R}^q)} \int \mathbf{E}^{\delta_x}(\|X_t - \eta_t(X_0)\|)\mu(dx)$$
$$= \sup_{x \in \mathbb{R}^q} \mathbf{E}^{\delta_x}(\|X_t - \eta_t(X_0)\|)$$
$$= \sup_{x \in \mathbb{R}^q} \mathbf{E}(\|\xi_t(x) - \eta_t(x)\|).$$

The proof is complete. □

We can now proceed with the proof of Theorem 3.9.

PROOF OF THEOREM 3.9. Let us fix two sequences $\{\mu_n\}, \{\nu_n\} \subset P(\mathbb{R}^q)$ so that $\|\mathbf{P}^{\mu_n}|_{\mathcal{F}^Y} - \mathbf{P}^{\nu_n}|_{\mathcal{F}^Y}\|_{\text{TV}} \to 0$, a constant $\alpha > 0$, and a function $f \in \text{Lip}(\mathbb{R}^q)$. In the following $\varepsilon_0, m, M > 0$ are as defined in Lemma 5.1, and we define

$$H_\varepsilon(x) := \frac{1}{\varepsilon}\int_0^\varepsilon h(\eta_s(x)) \, ds, \qquad \tilde{H}_\varepsilon := \frac{1}{\varepsilon}\int_0^\varepsilon h(X_s) \, ds.$$

By Lemma 5.2, we may choose $\varepsilon < \varepsilon_0$ such that

$$\sup_{s \leq \varepsilon} \sup_{\mu \in P(\mathbb{R}^q)} \mathbf{E}^\mu(\|X_s - \eta_s(X_0)\|) < \alpha.$$

By Lemma 5.1, there is an unbounded-Lipschitz function $g_\varepsilon$, with $\|g_\varepsilon\|_L \le m^{-1}$, such that $g_\varepsilon(H_\varepsilon(x)) = x$ for all $x \in \mathbb{R}^q$. In particular, we have for all $n \in \mathbb{N}$

$$\left|\int f\,d\mu_n - \int f\,d\nu_n\right| = |\mathbf{E}^{\mu_n}(f_\varepsilon(H_\varepsilon(X_0))) - \mathbf{E}^{\nu_n}(f_\varepsilon(H_\varepsilon(X_0)))|,$$

where we have written $f_\varepsilon := f \circ g_\varepsilon$. Now note that

$$\sup_{\mu \in P(\mathbb{R}^q)} |\mathbf{E}^\mu(f_\varepsilon(H_\varepsilon(X_0))) - \mathbf{E}^\mu(f_\varepsilon(\tilde{H}_\varepsilon))|$$

$$\le \|f_\varepsilon\|_L \|h\|_L \frac{1}{\varepsilon} \int_0^\varepsilon \sup_{\mu \in P(\mathbb{R}^q)} \mathbf{E}^\mu(\|X_s - \eta_s(X_0)\|)\,ds < \|h\|_L m^{-1} \alpha.$$

Therefore, we have for all $n \in \mathbb{N}$

$$\left|\int f\,d\mu_n - \int f\,d\nu_n\right| \le 2\|h\|_L m^{-1}\alpha + |\mathbf{E}^{\mu_n}(f_\varepsilon(\tilde{H}_\varepsilon)) - \mathbf{E}^{\nu_n}(f_\varepsilon(\tilde{H}_\varepsilon))|.$$

To proceed, note that

$$\frac{Y_\varepsilon}{\varepsilon} = \tilde{H}_\varepsilon + \frac{DV_\varepsilon}{\varepsilon}.$$

As $DV_\varepsilon/\varepsilon$ is Gaussian, its characteristic function vanishes nowhere. Therefore, using that $f_\varepsilon \in U_b(\mathbb{R}^q)$ and Proposition C.1, we may choose $u_\varepsilon \in U_b(\mathbb{R}^q)$ so that

$$\sup_{\mu \in P(\mathbb{R}^q)} |\mathbf{E}^\mu(f_\varepsilon(\tilde{H}_\varepsilon)) - \mathbf{E}^\mu(u_\varepsilon(Y_\varepsilon/\varepsilon))| < \alpha.$$

We thus obtain for every $n \in \mathbb{N}$

$$\left|\int f\,d\mu_n - \int f\,d\nu_n\right| \le 2\alpha(1 + \|h\|_L m^{-1}) + |\mathbf{E}^{\mu_n}(u_\varepsilon(Y_\varepsilon/\varepsilon)) - \mathbf{E}^{\nu_n}(u_\varepsilon(Y_\varepsilon/\varepsilon))|.$$

But as $\|\mathbf{P}^{\mu_n}|_{\mathcal{F}^Y} - \mathbf{P}^{\nu_n}|_{\mathcal{F}^Y}\|_{\mathrm{TV}} \to 0$, we evidently have

$$\limsup_{n \to \infty} \left|\int f\,d\mu_n - \int f\,d\nu_n\right| \le 2\alpha(1 + \|h\|_L m^{-1}).$$

Now note that $\alpha > 0$ and $f \in \mathrm{Lip}(\mathbb{R}^q)$ were arbitrary, so evidently

$$\left|\int f\,d\mu_n - \int f\,d\nu_n\right| \xrightarrow{n \to \infty} 0 \qquad \text{for all } f \in \mathrm{Lip}(\mathbb{R}^q).$$

The result follows from Corollary B.4. □

## APPENDIX A: MEASURABILITY OF PROBABILITY DISTANCES

The goal of this appendix is to prove that the distance $\|\mu - \nu\|_G$ between two probability kernels $\mu, \nu$ is measurable, provided that the family $G \subset C_b(S)$ is chosen appropriately. To this end we prove the following lemma.

LEMMA A.1. *Let $G \subset C_b(S)$ be uniformly bounded and equicontinuous. Then there is a countable collection $\{g_n : n \in \mathbb{N}\} \subset G$ such that*

$$\|\mu - \nu\|_G = \sup_n \left| \int g_n \, d\mu - \int g_n \, d\nu \right| \qquad \text{for all } \mu, \nu \in P(S).$$

PROOF. $P(S)$ is Polish, so there is a countable dense subset $\{\mu_n : n \in \mathbb{N}\} \subset P(S)$. As any probability measure on a Polish space is tight, we can find for every $n, m, p \in \mathbb{N}$ a compact set $K_{n,m,p} \subset S$ such that $\mu_n(K_{n,m,p}) > 1 - 1/p$ and $\mu_m(K_{n,m,p}) > 1 - 1/p$. Let us write $G_{n,m,p} = \{f|_{K_{n,m,p}} : f \in G\} \subset C_b(K_{n,m,p})$. By the Arzelà–Ascoli theorem the family $G_{n,m,p}$ is compact, and thus a forteriori separable, in the topology of uniform convergence. Therefore, we can find for every $n, m, p \in \mathbb{N}$ a countable family $\{g_k^{n,m,p} : k \in \mathbb{N}\} \subset G$ such that

$$\forall n, m, p \in \mathbb{N}, g \in G, \varepsilon > 0, \exists k \in \mathbb{N} \text{ s.t.} \sup_{x \in K_{n,m,p}} |g(x) - g_k^{n,m,p}(x)| < \varepsilon.$$

We claim that the countable family $G' = \{g_k^{n,m,p} : n, m, p, k \in \mathbb{N}\} \subset G$ satisfies

$$\|\mu - \nu\|_G = \sup_{g \in G'} \left| \int g \, d\mu - \int g \, d\nu \right| := \|\mu - \nu\|_{G'} \qquad \text{for all } \mu, \nu \in P(S).$$

Of course, the inequality $\|\mu - \nu\|_{G'} \leq \|\mu - \nu\|_G$ is trivial as $G' \subset G$, so it suffices to prove that for every $\mu, \nu \in P(S)$ there is a sequence $\{h_\ell : \ell \in \mathbb{N}\} \subset G'$ such that $|\mu(h_\ell) - \nu(h_\ell)| \to \|\mu - \nu\|_G$. To this end, let us fix $\mu, \nu \in P(S)$, and choose a sequence $\{h'_\ell : \ell \in \mathbb{N}\} \subset G$ such that $|\mu(h'_\ell) - \nu(h'_\ell)| \to \|\mu - \nu\|_G$. Note that

$$||\mu(h'_\ell) - \nu(h'_\ell)| - |\mu(h_\ell) - \nu(h_\ell)|| \leq |\mu(h'_\ell) - \mu(h_\ell)| + |\nu(h'_\ell) - \nu(h_\ell)|$$

by the reverse triangle inequality, so it suffices to find a sequence $\{h_\ell : \ell \in \mathbb{N}\} \subset G'$ such that $|\mu(h'_\ell) - \mu(h_\ell)| \to 0$ and $|\nu(h'_\ell) - \nu(h_\ell)| \to 0$. Fix $\ell \in \mathbb{N}$. By [23], Theorem II.6.8, we can choose $n, m \in \mathbb{N}$ such that $\|\mu_n - \mu\|_G < 1/\ell$ and $\|\mu_m - \nu\|_G < 1/\ell$. Choose $k \in \mathbb{N}$ such that $\sup_{x \in K_{n,m,\ell}} |h'_\ell(x) - g_k^{n,m,\ell}(x)| < 1/\ell$, and set $h_\ell = g_k^{n,m,\ell}$. Then we can estimate as follows:

$$|\mu(h'_\ell) - \mu(h_\ell)| \leq |\mu(h'_\ell) - \mu_n(h'_\ell)| + |\mu_n(h'_\ell) - \mu_n(h_\ell)|$$
$$\leq \|\mu_n - \mu\|_G + \left| \int_{K_{n,m,\ell}} (h'_\ell - h_\ell) \, d\mu_n \right| + \left| \int_{K^c_{n,m,\ell}} (h'_\ell - h_\ell) \, d\mu_n \right|$$

$$\leq \|\mu_n - \mu\|_G + \sup_{x \in K_{n,m,\ell}} |h'_\ell(x) - h_\ell(x)|$$
$$+ 2\sup_{g \in G} \|g\|_\infty \mu_n(K^c_{n,m,\ell})$$
$$\leq \frac{2}{\ell}\Big(1 + \sup_{g \in G} \|g\|_\infty\Big),$$

where $K^c$ denotes the complement of a set $K$. The identical bound is found for $|\nu(h'_\ell) - \nu(h_\ell)|$. Repeating the procedure for every $\ell \in \mathbb{N}$, we evidently construct a sequence $\{h_\ell\}$ with the desired properties. This completes the proof. □

This result will be used in the following fashion.

COROLLARY A.2. *Let $(\Omega, \mathcal{F})$ be a measurable space and let $\mu:\Omega \times \mathcal{B}(S) \to [0,1]$ and $\nu:\Omega \times \mathcal{B}(S) \to [0,1]$ be probability kernels. Moreover, let $G \subset C_b(S)$ be uniformly bounded and equicontinuous. Then $\|\mu - \nu\|_G$ is a random variable [i.e., the map $\omega \mapsto \|\mu(\omega, \cdot) - \nu(\omega, \cdot)\|_G$ is measurable].*

PROOF. Immediate from the previous lemma. □

Corollary A.2 is used implicitly throughout the paper without further comment.

## APPENDIX B: MERGING OF PROBABILITY MEASURES

It is well known that a sequence of probability measures $\{\mu_n\} \subset P(S)$ converges weakly to $\mu \in P(S)$ if and only if $\|\mu_n - \mu\|_{\text{BL}} \to 0$ [13], Theorem 11.3.3. In particular, as the class of bounded-Lipschitz functions is uniformly dense in $U_b(S)$, it follows that if $\mu_n(f) \to \mu(f)$ for all $f \in \text{Lip}(S)$ then $\|\mu_n - \mu\|_{\text{BL}} \to 0$. This is in some sense surprising; evidently the convergence of the expectation of every function $f \in \text{Lip}(S)$ separately already implies that this convergence holds uniformly over $\text{Lip}(S)$, without any further assumptions.

The purpose of this appendix is to show that a similar statement holds for the *merging* of two sequences of probability measures. This result was already proved in [22] and in [7], Section 6, for probability measures on any Polish space. We provide here an alternative and much simpler proof, which is however restricted to probability measures on $\mathbb{R}^d$. In this paper only the latter will be needed.

PROPOSITION B.1. *Let $\{\mu_n\}, \{\nu_n\} \subset P(\mathbb{R}^d)$ satisfy*
$$\left|\int f \, d\mu_n - \int f \, d\nu_n\right| \xrightarrow{n \to \infty} 0 \quad \text{for all } f \in \text{Lip}(\mathbb{R}^d).$$

Moreover, let $G \subset U_b(\mathbb{R}^d)$ be a uniformly bounded and uniformly equicontinuous family of functions. Then $\|\mu_n - \nu_n\|_G \xrightarrow{n \to \infty} 0$.

The proof is based on the following well-known result from Banach space theory, which states in essence that the claim is true for probability measures on $\mathbb{N}$ (rather than $\mathbb{R}^d$). An elementary proof can be found in [1], Theorem 4.32.

LEMMA B.2 (Schur property of $\ell_1$). *A sequence in $\ell_1$ converges in the weak topology if and only if it converges in the norm topology.*

REMARK B.3. Note that this result only holds for *sequences*. Indeed, it can not hold for nets, as that would imply that the weak and norm topologies coincide.

We now turn to the proof of Proposition B.1. The basic idea is to reduce to the setting of Lemma B.2 by introducing a partition of unity.

PROOF OF PROPOSITION B.1. For every $\alpha > 0$, define the countable family of functions $V^\alpha = \{x \mapsto \varphi_{k_1}(\alpha x^1) \cdots \varphi_{k_d}(\alpha x^d) : (k_1, \ldots, k_d) \in \mathbb{Z}^d\} \subset U_b(\mathbb{R}^d)$, where $\varphi_k(x) = \cos^2(\pi(x-k)/2) I_{|x-k| \leq 1}$. The following facts are easily verified:

1. $0 \leq \varphi(x) \leq 1$ for all $\varphi \in V^\alpha$;
2. For every $x \in \mathbb{R}$, at most $N$ elements of $\varphi \in V^\alpha$ satisfy $\varphi(x) > 0$ (where $N \in \mathbb{N}$ depends only on the state space dimension $d$);
3. For every $x \in \mathbb{R}^d$, we have $\sum_{\varphi \in V^\alpha} \varphi(x) = 1$;
4. $\sup_{\varphi \in V^\alpha} \|\varphi\|_L < \infty$ (i.e., $V^\alpha$ is equilipschitzian).

Thus $V^\alpha$ is a partition of unity of $\mathbb{R}^d$ with some additional uniformity properties (which will be important in the following).

Fix $\varepsilon > 0$. As $G$ is uniformly equicontinuous, there exists a $\delta > 0$ such that $\|x - y\| \leq \delta$ implies $|g(x) - g(y)| \leq \varepsilon$ for all $g \in G$. Choose $\alpha$ large enough so that every element of $V^\alpha$ is supported inside a ball of radius $\delta$. Moreover, choose for every $\varphi \in V^\alpha$ an arbitrary point $x_\varphi \in \mathbb{R}^d$ in the support of $\varphi$. Then

$$\|\mu_n - \nu_n\|_G = \sup_{g \in G} \left| \sum_{\varphi \in V^\alpha} \int g\varphi \, d\mu_n - \sum_{\varphi \in V^\alpha} \int g\varphi \, d\nu_n \right|$$

$$\leq \sup_{g \in G} \sum_{\varphi \in V^\alpha} \left| \int g\varphi \, d\mu_n - \int g\varphi \, d\nu_n \right|$$

$$\leq \sup_{g \in G} \sum_{\varphi \in V^\alpha} \left\{ \left| \int (g - g(x_\varphi))\varphi \, d\mu_n \right| + \left| \int (g - g(x_\varphi))\varphi \, d\nu_n \right| \right\}$$

$$+ \sup_{g \in G} \sum_{\varphi \in V^\alpha} \left\{ |g(x_\varphi)| \left| \int \varphi \, d\mu_n - \int \varphi \, d\nu_n \right| \right\}$$

$$\leq 2\varepsilon + \sup_{g \in G} \|g\|_\infty \sum_{\varphi \in V^\alpha} \left| \int \varphi \, d\mu_n - \int \varphi \, d\nu_n \right|.$$

Suppose that we can show that for any $\alpha > 0$

$$\sum_{\varphi \in V^\alpha} \left| \int \varphi \, d\mu_n - \int \varphi \, d\nu_n \right| \xrightarrow{n \to \infty} 0.$$

Then the result follows as $\varepsilon > 0$ was arbitrary.

Now note that $V^\alpha$ is countable, so it may be ordered as $V^\alpha = \{\chi_k : k \in \mathbb{N}\}$. For any finite signed measure $\rho$, the sequence $(\rho(\chi_k))_{k \in \mathbb{N}} \in \ell_1$. We must establish that $(\mu_n(\chi_k) - \nu_n(\chi_k))_{k \in \mathbb{N}}$ converges to zero in the $\ell_1$-norm. Therefore, by Lemma B.2, it suffices to prove that this convergence holds in the weak topology. In particular, define for every $z \in \ell_\infty$ the function $f_z := \sum_k z_k \chi_k$. Then it suffices to show that

$$\left| \int f_z \, d\mu_n - \int f_z \, d\nu_n \right| \xrightarrow{n \to \infty} 0 \qquad \text{for every } z \in \ell_\infty.$$

By our assumptions this is the case if $\|f_z\|_\infty < \infty$ and $\|f_z\|_L < \infty$. But this holds for any $z \in \ell_\infty$; indeed, it is easily seen from the properties of $V^\alpha$ that $\|f_z\|_\infty \leq \|z\|_\infty$ and $\|f_z\|_L \leq 2N\|z\|_\infty \sup_{\varphi \in V^\alpha} \|\varphi\|_L$. This completes the proof. □

The following corollary is immediate.

COROLLARY B.4. *For $\{\mu_n\}, \{\nu_n\} \subset P(\mathbb{R}^d)$, the following are equivalent:*

1. $|\int f \, d\mu_n - \int f \, d\nu_n| \to 0$ *as* $n \to \infty$ *for every* $f \in \mathrm{Lip}(\mathbb{R}^d)$;
2. $\|\mu_n - \nu_n\|_{\mathrm{BL}} \to 0$ *as* $n \to \infty$.

It should be noted that for sequences of probability measures in $P(K)$, where $K$ is a compact Polish space, the same results can be proved in a completely elementary fashion; indeed, in this case any uniformly bounded and equicontinuous family $G \subset C_b(K)$ is compact in the topology of uniform convergence (by the Arzelà–Ascoli theorem), so that it can be covered by a finite number of arbitrarily small balls. The previous results then follow from elementary arguments. When the state space is not compact, however, the result is far from obvious and relies heavily on the (elementary but nontrivial) Schur property of $\ell_1$.

## APPENDIX C: UNIFORM APPROXIMATION AND CONVOLUTION

In order to verify uniform observability for additive noise observation models, the following result is of central importance. Though the result seems to be of independent interest—it states that the range of a convolution operator on $U_b(\mathbb{R}^d)$ is uniformly dense in $U_b(\mathbb{R}^d)$ under a mild condition—the author was not able to find a statement or proof of this result in the literature.

PROPOSITION C.1. *Suppose that the characteristic function of the probability measure $\mu \in P(\mathbb{R}^d)$ vanishes nowhere. Then the family $\{f * \mu : f \in U_b(\mathbb{R}^d)\} \subset U_b(\mathbb{R}^d)$ is dense in $U_b(\mathbb{R}^d)$ in the topology of uniform convergence.*

The difficulty here is that we seek *uniform* density of functions on a noncompact space; as the Banach space dual $U_b(\mathbb{R}^d)^*$ contains elements which are not countably additive, this precludes the routine application of the Hahn–Banach theorem (see [14] for related results). We circumvent this problem by using the elementary properties of convolutions to "push" the approximation problem into $L^1(\mathbb{R}^d)$, where standard approximation results are readily available.

PROOF OF PROPOSITION C.1. We first collect some well-known facts about convolutions. Let $\varphi \in L^1(\mathbb{R}^d)$ be any function such that $\int \varphi(x)\, dx = 1$. Define $\varphi^t(x) = t^{-d}\varphi(t^{-1}x)$. Then by [16], Theorem 8.14, we have $\|f * \varphi^t - f\|_\infty \to 0$ as $t \to 0$ for any $f \in U_b(\mathbb{R}^d)$. Moreover, for $f \in U_b(\mathbb{R}^d)$ and $\varrho \in L^1(\mathbb{R}^d)$, we have by [16], Proposition 8.8, that $\|f * \varrho\|_\infty \leq \|f\|_\infty \|\varrho\|_1$ and $f * \varrho \in U_b(\mathbb{R}^d)$. Finally, if $\varrho \in L^1(\mathbb{R}^d)$ and $\nu \in P(\mathbb{R}^d)$, then $\varrho * \nu \in L^1(\mathbb{R}^d)$ by [16], Proposition 8.49.

Fix $\varphi \in L^1(\mathbb{R}^d)$ as above and let $f \in U_b(\mathbb{R}^d)$ and $k \in \mathbb{N}$. Then we may choose $t > 0$ such that $\|f * \varphi^t - f\|_\infty \leq k^{-1}$. Now suppose that we can find a function $\varrho_k \in L^1(\mathbb{R}^d)$ such that $\|\varphi^t - \varrho_k * \mu\|_1 \leq k^{-1}$. Then we can evidently estimate

$$\|f - f * (\varrho_k * \mu)\|_\infty \leq \|f - f * \varphi^t\|_\infty + \|f\|_\infty \|\varphi^t - \varrho_k * \mu\|_1 \leq k^{-1}(1 + \|f\|_\infty).$$

But note that $f * (\varrho_k * \mu) = (f * \varrho_k) * \mu$ and $g_k := f * \varrho_k \in U_b(\mathbb{R}^d)$. Repeating the procedure for every $k \in \mathbb{N}$, we find a sequence $\{g_k\} \subset U_b(\mathbb{R}^d)$ such that $\|f - g_k * \mu\|_\infty \to 0$. As the function $f \in U_b(\mathbb{R}^d)$ was arbitrary, the result follows.

It thus remains to show that for every $t > 0$ and $k \in \mathbb{N}$, we can find a function $\varrho \in L^1(\mathbb{R}^d)$ such that $\|\varphi^t - \varrho * \mu\|_1 \leq k^{-1}$. It suffices to show that the family $\{\varrho * \mu : \varrho \in L^1(\mathbb{R}^d)\} \subset L^1(\mathbb{R}^d)$ is dense in $L^1(\mathbb{R}^d)$. To this end, consider

$$T\Phi := \operatorname{span}\{x \mapsto \Phi(x - a) : a \in \mathbb{R}^d\} \subset L^1(\mathbb{R}^d), \qquad \Phi(x) := e^{-\|x\|^2/2}.$$

Evidently $\{\varrho * \mu : \varrho \in T\Phi\}$ is the span of all translates of the function $\Phi * \mu$. But the Fourier transform $(\Phi * \mu)^\wedge = \Phi^\wedge \mu^\wedge$ vanishes nowhere, so $\{\varrho * \mu : \varrho \in T\Phi\}$ is dense in $L^1(\mathbb{R}^d)$ by the Wiener Tauberian theorem [25], Theorem 9.5. $\square$

As an application, we prove the following result.

PROPOSITION C.2. *Let $\{\mu_n\}, \{\nu_n\} \subset P(\mathbb{R}^d)$, and let $\xi \in P(\mathbb{R}^d)$ be a probability measure whose characteristic function vanishes nowhere. Then*

$$\|\mu_n * \xi - \nu_n * \xi\|_{\mathrm{BL}} \xrightarrow{n \to \infty} 0 \quad \text{if and only if} \quad \|\mu_n - \nu_n\|_{\mathrm{BL}} \xrightarrow{n \to \infty} 0.$$

In other words, if $\xi$ is a probability measure whose characteristic function vanishes nowhere, then the convolution operator $C_\xi : P(\mathbb{R}^d) \to P(\mathbb{R}^d)$ defined as $C_\xi \mu = \mu * \xi$ is uniformly continuous, injective and the inverse operator $C_\xi^{-1} : \mathrm{Ran}\, C_\xi \to P(\mathbb{R}^d)$ is uniformly continuous (relative to the $\|\cdot\|_{\mathrm{BL}}$-norm).

PROOF OF PROPOSITION C.2. Denote by $\bar\xi$ the reflected probability measure defined by

$$\int f(x)\bar\xi(dx) = \int f(-x)\xi(dx) \qquad \text{for all } f \in B(\mathbb{R}^d).$$

Clearly the characteristic function of $\bar\xi$ vanishes nowhere. Now note that we obtain for any probability measure $\mu \in P(\mathbb{R}^d)$ the identity

$$\int f(x)(\mu * \xi)(dx) = \int (f * \bar\xi)(x)\mu(dx) \qquad \text{for all } f \in B(\mathbb{R}^d).$$

Moreover, it is easily verified that $f * \bar\xi \in U_b(\mathbb{R}^d)$ whenever $f \in U_b(\mathbb{R}^d)$.

Let us first suppose that $\|\mu_n - \nu_n\|_{\mathrm{BL}} \to 0$. Then

$$\left| \int f\, d\mu_n - \int f\, d\nu_n \right| \xrightarrow{n \to \infty} 0 \qquad \text{for all } f \in U_b(\mathbb{R}^d)$$

as the family of bounded-Lipschitz functions is dense in $U_b(\mathbb{R}^d)$ in the topology of uniform convergence [12], Lemma 8. Therefore

$$\left| \int f\, d(\mu_n * \xi) - \int f\, d(\nu_n * \xi) \right| = \left| \int (f * \bar\xi)\, d\mu_n - \int (f * \bar\xi)\, d\nu_n \right| \xrightarrow{n \to \infty} 0$$

for every $f \in U_b(\mathbb{R}^d)$. That $\|\mu_n * \xi - \nu_n * \xi\|_{\mathrm{BL}} \to 0$ follows from Corollary B.4.

Conversely, let us suppose that $\|\mu_n * \xi - \nu_n * \xi\|_{\mathrm{BL}} \to 0$, so that

$$\left| \int f\, d\mu_n - \int f\, d\nu_n \right| \xrightarrow{n \to \infty} 0, \qquad \text{whenever } f \in \{g * \bar\xi : g \in U_b(\mathbb{R}^d)\}.$$

By Proposition C.1, the family $\{g * \bar\xi : g \in U_b(\mathbb{R}^d)\}$ is uniformly dense in $U_b(\mathbb{R}^d)$; therefore this convergence holds for any $f \in U_b(\mathbb{R}^d)$. But then Corollary B.4 implies that $\|\mu_n - \nu_n\|_{\mathrm{BL}} \to 0$, and the proof is complete. $\square$

DEPARTMENT OF OPERATIONS RESEARCH
AND FINANCIAL ENGINEERING
PRINCETON UNIVERSITY
PRINCETON, NEW JERSEY 08544
USA
E-MAIL: rvan@princeton.edu